\documentclass[a4paper,12pt]{amsart}
\usepackage{amssymb,amsfonts,amsmath}

\usepackage[T1]{fontenc} 
\usepackage[latin1]{inputenc} 

\def\d{\delta}
\def\H{\mathcal{H}}

\def\C{\mathbb{C}}
\def\c2{\mathbb{C}^2}
\def\R{\mathbb{R}}

\def\Z{\mathbb{Z}}
\def\N{\mathbb{N}}
\def\P{\mathbb{P}}

\def\1{\bold{1}}

\def\a{\alpha}
\def\b{\beta}
\def\e{\varepsilon}
\def\l{\lambda}
\def\f{\varphi}
\def\g{\gamma}
\def\p{\psi}
\def\r{\varrho}
\def\o{\omega}

\def\D{\overline{\partial}}
\def\dist{\operatorname{dist}}
\def\supp{\operatorname{supp}}
\def\H11{H^{1,1}_\R}
\def\I{\mathcal I}

\newcommand{\pair}[2]{\left\langle #1,#2 \right\rangle}
\newtheorem{lem}{Lemma}[section]
\newtheorem{pro}[lem]{Proposition}
\newtheorem{defi}[lem]{Definition}
\newtheorem{def/not}[lem]{Definition/Notations}
\newtheorem{Pro/def}[lem]{Proposition/Definition}

\newtheorem{thm}[lem]{Theorem}
\newtheorem{cor}[lem]{Corollary}
\newtheorem{rqe}[lem]{Remark}

\newtheorem{exa}[lem]{Example}
\newtheorem{exas}[lem]{Examples}

\newtheorem{ques}[lem]{Question}

\newenvironment{proof.3.1}
{\vskip .2cm \noindent {{\it Proof of theorem 3.1.}}}{\hfill $\Box$}
\newenvironment{proof.3.4}
{\vskip .2cm \noindent {{\it Proof of theorem 3.4.}}}{\hfill $\Box$}

\begin{document}

\title[Regularity of dynamical Green functions]
{Regularity of dynamical Green functions}

\author{Jeffrey DILLER \& Vincent GUEDJ}

\thanks{The first author gratefully acknowledges support from the National Science Foundation during the preparation of this article.}

\begin{abstract}
For meromorphic maps of complex manifolds, ergodic theory and pluripotential
theory are closely related.  In nice enough situations, dynamically defined 
Green's functions give rise to invariant currents which intersect to yield 
measures of maximal entropy.  `Nice enough' is often a condition
on the regularity of the Green's function.  In this paper we look at a 
variety of regularity properties that have been considered for dynamical 
Green's functions.  We simplify and extend some known results and prove
several others which are new.  We also give some examples indicating
the limits of what one can hope to achieve in complex dynamics by relying 
solely on the regularity of a dynamical Green's function.
\end{abstract}

\maketitle

\section*{Introduction}

A holomorphic, or more generally, meromorphic self-map $f:X\to X$ of a compact
complex 
manifold $X$ induces 
actions $f^*,f_*:H^*(X,\R)\to H^*(X,\R)$ on the real cohomology groups of
$X$.  It is conjectured that when these actions are suitably well-behaved,
then the topological entropy $h_{top}(f)$ of $f$ should be 
$\log \rho(f^*)$, where
$\rho(\cdot)$ denotes spectral radius.  This conjecture has motivated a great
deal of research in the past fifteen years, and it has been verified in
some important cases (see [Gr], [Y], [Sm], [FS 2], [Du], [G 3]).  
It is known, for instance, that
the inequality $h_{top}(f) \leq \log\rho(f)$ always holds [DS 1].

The main strategy for proving the reverse inequality has been to look for
an invariant measure whose metric entropy is maximal, i.e. equal to 
$\log\rho(f)$.
However, rather than try to realize the measure directly from the dynamics
of $f$, it often seems more promising to use the dynamics to construct 
invariant positive closed currents and then try to obtain the measure as an
intersection of these currents.  The drawback is that in passing from currents 
to measures, one must somehow make sense of what is essentially a product
of distributions.  For positive closed currents, this is usually done by
resorting to `potentials' for the currents and integrating by parts. 
Success depends on having potentials that are substantially more regular
than the currents themselves. 
The purpose of this paper is to better understand the regularity properties
of potentials associated to dynamically-defined positive closed $(1,1)$ 
currents.  Such potentials will be functions, the \emph{dynamical Green's 
functions} in the title of the paper.  

In the first section we describe the best possible situation: 
\emph{holomorphic} maps.  We present a simple proof, due to Dinh and Sibony 
(see [DS 2], Theorem 3.7.1; also [DS 3], Proposition 2.4), of the fact
that a dynamical Green's function associated to a holomorphic map must
be H\"older continuous with H\"older exponent controlled by what we 
call the
\emph{topological Lyapunov exponent}
$$
\chi_{top}(f):=\lim_{n \rightarrow +\infty} \frac{1}{n} \log \sup_{x \in \P^k}
||D_x f^n||
$$
of the map.  A straightforward example shows that this estimate is sharp.  

In the remaining sections of the paper, we turn to the more general class of 
meromorphic self-maps, a principal goal being to see what remains of 
H\"older continuity for the Green's function once one leaves the holomorphic 
setting.  

Our first result, proven in section 2, is a general existence theorem 
for dynamical Green's functions of meromorphic maps in any dimension.  If
$f:X\to X$ is a meromorphic map of a compact K\"ahler manifold, then we say
that $f$ is \emph{1-stable} if the induced action $f^*$ on $H^{1,1}_\R(X)$
satisfies $(f^n)^* = (f^*)^n$ for all $n\in\N$.  Given a class 
$\eta\in H^{1,1}_\R(X)$ satisfying $f^*\eta = \lambda\eta$ and a smooth form
$\omega$ representing $\eta$, one can try to construct an invariant current
representing $\eta$ as follows.  For each $n\in\N$ we have an $L^1$ function 
$g^\omega_n:X\to \R\cup\{\pm\infty\}$ satisfying
$$
\lambda^{-n}f^{n*}\omega - \omega = dd^c g^\omega_n.
$$
If $g^\omega_n$ converges in $L^1$ to some function $g^\omega$ (the dynamical
Green's function), then the current $T_\eta := \omega + dd^c g^\omega$ 
automatically satisfies $f^* T = \lambda T$.

\vskip.3cm
\noindent {\bf Theorem A.}
{\it 
Suppose that $f:X\to X$ is a 1-stable meromorphic map of a compact K\"ahler 
surface and that the induced action $f^*$ has a unique simple
eigenvalue $\lambda$ of largest modulus with eigenspace generated by a nef
class $\eta$.  Then for any smooth form $\omega$ representing $\eta$, it can
be arranged that the sequence $(g^\omega_n)$ is decreasing and $L^1$
convergent.  The closed current $T := \omega + dd^c g^\omega$ is positive
and independent of $\omega$.}
\vskip.3cm

Our proof 
follows Sibony [S] who considered the case $X=\P^k$ and Guedj [G 1,2] who 
considered general $X$.
The novelty here is that we do not assume that the smooth representative 
$\omega$ can be chosen positive.  Hence it requires some new ideas to
establish that the sequence approximating $g^\omega$ is decreasing and to
show that the invariant current $T$ is positive.  We remark that in dimension
two, the theorem applies to nearly all reasonable meromorphic maps
(see Corollary \ref{dim2cor}).  After section 2, we restrict attention to
maps of complex \emph{surfaces}.

When the map $f$ in Theorem A is not holomorphic, the Green's 
function $g_\omega$ will not be continuous.  
It will typically have a logarithmic pole at each point of 
indeterminacy for $f$ and its iterates.  We let $\mathcal{I}_f$ denote the
closure of the set of all such points.  Though $\mathcal{I}_f$ can be all of 
$X$, as the first example in Section 6 shows, there are many situations
where the complement of $\mathcal{I}_f$ is large, and one can then hope for
continuity of $g^+$ in $X- \mathcal{I}_f$.  In section 3, we validate this
hope in some interesting special cases.  Indeed, we give a unified approach
to proving something analogous to, but weaker than, H\"older continuity for
$g^+$ for some large classes of birational surface maps (Theorem 3.1 and the
comment following its proof) and of polynomial maps of $\C^2$ (Theorem 3.4).  
We point out concerning Theorem 3.4 that in the important case where the first 
dynamical degree $\lambda$ exceeds the topological degree, we know of no
examples where the hypothesis of the theorem fails.

In Section 4, we consider a still weaker regularity condition for birational 
surface maps.  We state a quantitative recurrence property for 
points of indeterminacy that turns out to be equivalent to the condition that 
the derivative $dg^+$ of the Green's function be in $L^2$.  A similar, 
slightly stronger
$L^2$ condition has been used with much success in [BD] and [Du] to produce
measures of maximal entropy for birational maps.  With our version, the
construction of the measure still succeeds but its fine dynamical properties 
remain unclear; in particular, we do not know if $\log||Df||$ is integrable
with respect to the measure, a property that is important for applying
many of the theorems and techniques from smooth ergodic theory.

Continuing with birational surface maps in Section 5, we consider 
what is perhaps the weakest relevant regularity 
condition of all: $g^+$ is integrable with respect to (the trace measure of)
$T^-$, the invariant current associated to $f^{-1}$.  This condition
guarantees that $\mu = T^+\wedge T^-$ is a well-defined probability 
measure.  Indeed with no further assumption on $f$ (i.e. on $g^+$), we
prove the following.

\vskip.3cm
\noindent {\bf Theorem B.}
{\it The measure $\mu_f$ is $f$-invariant and mixing, and it does not charge any compact complex curve.}
\vskip.3cm

The proof that $\mu_f$ does not charge curves is distinctly indirect, 
depending on among other things, a characterization (Proposition 
\ref{finiteg}) of the $L^2$ condition used in [BD].

We present several telling examples throughout the paper, and Section 6 is
devoted to two of these.  The first shows that indeterminacy orbits of a 
meromorphic map can be dense.  That is, $\mathcal I_f = X$.  
The second builds on an example
due to Favre [F] and demonstrates that one can
have $g^+$ integrable with respect to $T^-$ without necessarily having that
$dg^+$ is in $L^2$.  

\section{Holomorphic maps}

Let $f:\P^k \rightarrow \P^k$ be a holomorphic endomorphism of the complex 
projective space $\P^k$. In homogeneous coordinates
$f=[P_0:\ldots:P_k]$ where the $P_j$'s are homogeneous polynomials of the same
degree ${\bf \l}$ with no common zero outside the origin. 
We shall always assume $\l \geq 2$.

Let $\o$ denote the Fubini-Study K\"ahler form on $\P^k$. Then $f^* \o$ is a well-defined
smooth closed $(1,1)$-form on $\P^k$ which is cohomologous to $\l \o$.  Thus
it follows from the $dd^c$-lemma (see [GH] p149) that
$$
\frac{1}{\l}f^*\o=\o+dd^c \g,
$$
where $\g \in {\mathcal C}^{\infty}(\P^k)$ is uniquely determined up to an 
additive constant.  Here $d=\partial+\D$ and 
$d^c=\frac{1}{2i\pi}(\D-\partial)$.  Pulling back the previous equation by 
$f^n$ yields 
$$
\frac{1}{\l^n}(f^n)^* \o=\o+dd^c g_n, \text{ where }
g_n=\sum_{j=0}^{n-1} \frac{1}{\l^j} \g \circ f^j.
$$
The sequence of positive closed $(1,1)$-forms
$\l^{-n}(f^n)^*\o$ converges weakly to the so-called Green current
$$
T_f=\o+dd^c g_f, \text{ where } g_f:=\sum_{j \geq 0} \frac{1}{\l^j} \g \circ f^j.
$$
This is a dynamically interesting current. It was constructed by H.Brolin [Bro] 
(polynomial case) and M.Lyubich [Ly] (rational case) when $k=1$, and by 
Hubbard-Papadopol [HP] and Fornaess-Sibony [FS 2] in higher dimensions. We
refer the reader to [S]
for its basic properties. Our aim here is to give a very simple proof of the fact
that the (dynamical) Green function $g_f$ is H\"older continuous. To this end we introduce
the {topological Lyapunov exponent} of $f$,
$$
\chi_{top}(f):=\lim_{n \rightarrow +\infty} \frac{1}{n} \log \sup_{x \in \P^k}
||D_x f^n||.
$$
That the limit exists follows from the submultiplicativity of the sequence
$(\sup_{x \in \P^k} ||D_xf^n||)$.
The definition clearly does not depend on the choice of the norm $|| \cdot ||$.
Also, the supremum in the definition can be considered only on the Julia set
$J_f$ of $f$. Let us recall that the {\it Fatou set} ${\mathcal F}_f$ of $f$ is the
largest open subset of $\P^k$ on which the sequence of iterates $(f^n)$ forms
a normal family. The {\it Julia set} $J_f$ is the complement of the Fatou set.

The next result, together with its proof, is essentially Theorem 3.7.1
in [DS 2] (see Proposition 2.4 in [DS 3] for a more general
statement).  We include it here for the convenience of the
reader, because many of the results in the following sections may be viewed
as attempts to salvage what remains of Theorem 1.1 when one passes from 
holomorphic to more badly behaved meromorphic maps.

\begin{thm}
The Green function $g_f$ is H\"older continuous of exponent $\a>0$, for every
$\a < \log \l/\chi_{top}(f)$.
\end{thm}

\begin{proof}
Set $M=\sup_{x \in \P^k} ||D_x f||$. A straightforward induction yields,
for all $x,y \in \P^k$ and all $j \in \N$,
$$
d(f^j x, f^j y) \leq M^j d(x,y).
$$
Here $d$ denotes the distance associated to the Fubini-Study metric on $\P^k$.
Since $\g$ is smooth, it is in particular H\"older-continuous of exponent 
$\a>0$, for any $\a \leq 1$.  We fix $\a<\log \l/M$ and estimate
$$
|g_f(x)-g_f(y)| \leq \sum_{j \geq 0} \frac{1}{\l^j} | \g \circ f^j(x)-\g \circ f^j(y)|
\leq C_{\a} d(x,y)^{\a},
$$
where $C_{\a}=\sum_{j \geq 0} \l^{-j} M^{\a j}<+\infty$.

Replacing $f$ by $f^n$ in the above argument lowers the constant $M$ to
$M_n=$ $(\sup ||D_x f^n||)^{1/n}$. 
Letting $n \rightarrow +\infty$ yields the desired upper-bound.
\end{proof}
\vskip.1cm

Example 1.2 shows that the bound in this theorem is optimal.
One can also establish bounds in the other direction using the infimum of the 
differential on the Julia set.  These imply in particular that the affine Green functions $G_c$ of
quadratic maps $f_c(z)=z^2+c$ with $c \in \R$ are H\"older continuous of 
exponent $\a_c$ with $\a_c \rightarrow 0$ as $c \rightarrow +\infty$.

\begin{exas}
Consider the quadratic family of holomorphic endomorphisms of the Riemann sphere
$f_c:\P^1 \rightarrow \P^1$, given by quadratic polynomials in some affine chart,
$$
f_c(z)=z^2+c.
$$
We let $G_c(z)=g_{f_c}[1:z]$ denote the affine Green function.

 1) If $c=0$ one easily computes $G_0(z)=\log^+|z|$,
$J_{f_0}=S^1$ is the unit circle, and $\chi_{top}(f_0)=\log2$.
If $c=-2$, then $f_{-2}$ is semi-conjugate to $f_0$, and one can compute
its iterates explicitly. This yields $J_{f_{-2}}=[-2,2]$,
$\chi_{top}(f_{-2})=2\log 2$ and
$$
g_{-2}(z)=\log \max \left( \left|\frac{z+\sqrt{z^2-4}}{2}\right| ;
\left|\frac{z-\sqrt{z^2-4}}{2}\right| \right).
$$
Observe that this is a H\"older-continuous function of exponent $1/2$.

 2) If $|c| \leq 2$, one can easily show that the Julia set
$J_{f_c}$ is always contained in the closed disk centered at the origin and of 
radius $2$.  We infer $\chi_{top}(f_c) \leq 2 \log 2$. Note that this disk 
contains the Mandelbrot set ${\mathcal M}$, i.e. the set of parameter values 
$c$ such that $J_{f_c}$ is connected.

More generally, if $f$ is any monic centered polynomial
of degree $\l$ with connected Julia set, it was proved by X.Buff [Bu] that 
$$
\chi_{top}(f) \leq 
\sup_{J_f} |f'| \leq 2 \log \l.
$$
We also have bounds from below. If $\mu$ is any invariant ergodic measure 
such that $\log ||(Df)^{\pm 1}|| \in L^1(\mu)$, then its Lyapunov exponent
$\chi_{\mu}(f)$ satisfies $\chi_{\mu}(f) \leq \chi_{top}(f)$.
In particular if $\mu=\o+dd^c g_f=dd^c G_f$ is the Brolin-Lyubich measure, then
$$
\log\l \leq \log \l+\sum_{f'(c)=0} G_f(c)=\chi_{\mu}(f) \leq \chi_{top}(f).
$$
\end{exas}

\begin{rqe}
For $X=\P^1$, H\"older continuity of dynamical Green functions was 
first established by N.Sibony (see [CG] Theorem 8.3.2). It was then 
generalized to endomorphisms of $\P^k$ by J.-Y.Briend [Bri] and M.Kosek [K]. 

As we explain below, the elementary proof given 
above applies to other manifolds.  Slightly modifying the proof shows also that 
if $(f_t)_{t \in M}$ is a holomorphic family of endomorphisms of the 
same degree $\l$, then the Green function $(x,t) \mapsto g_{f_t}(x)$ is 
H\"older continuous with respect to the parameter $t$.
\end{rqe}

Consider a holomorphic endomorphism $f:X \rightarrow X$ of some compact 
K\"ahler manifold $X$.  Then $f^*$ respects complex conjugation and bidegree 
of forms, and therefore restricts to a linear action on 
$H_\R^{1,1}(X) := H^2(X,\R)\cap H^{1,1}(X)$. 

Assume that $f^*\eta=\l \eta$ for some $\l>1$ and $\eta\in\H11(X)$. Then
if $\o$ is a smooth closed $(1,1)$ form representing $\eta$, we have 
$\g\in{\mathcal C}^{\infty}(X,\R)$ such that 
$$
\frac{1}{\l}f^* \o=\o+dd^c \g.
$$ 
Pulling back by $f^n$ yields
$$
\frac{1}{\l^n} (f^n)^* \o \longrightarrow T_{\eta}=\o+dd^c g_{\o},
\text{ where }
g_{\o}=\sum_{j \geq 0} \frac{1}{\l^j} \g \circ f^j.
$$
Observe that the dynamical Green current $T_{\eta}$ only depends on $\eta$: if
$\o'$ also represents $\eta$, then 
$\o'=\o+dd^c u$ for some smooth function $u$. Hence
$u \circ f^n/\l^n\to 0$ uniformly on $X$, and
  $$ \frac{1}{\l^n} (f^n)^* \o'
  =\frac{1}{\l^n} (f^n)^* \o+
  dd^c \left( \frac{1}{\l^n} u \circ f^n \right) \rightarrow T_{\eta}.
  $$
The same proof as above shows that
{\bf the Green function  $g_{\o}$ is H\"older continuous}.

In concluding this section, we recall (see e.g. [Meo], section I) that
if $S$ is any positive closed $(1,1)$ current on $X$, then we may
define the pullback $f^* S$.  Namely, we write $S = \eta + dd^c u$ where
$\eta$ is smooth closed $(1,1)$ form and $u$ is a quasi-plurisubharmonic
(henceforth `qpsh') function, and we set $f^* S = f^*\eta + dd^c
u\circ f$.  The result is another positive closed $(1,1)$ current on
$X$.  

A cohomology class $\eta$ is \emph{pseudoeffective} if it can be represented
by a {\it positive} closed current $S$ of bidegree $(1,1)$.  In this case, 
$f^*S$ is a well-defined positive closed current of bidegree $(1,1)$ on $X$ 
which represents $f^*\eta$.  Thus $f^*$ preserves the cone
$H^{1,1}_{psef}(X)\subset\H11(X)$ of pseudoeffective classes. Because the 
pseudoeffective cone is closed, convex and strict (i.e. contains no lines), it follows 
from Perron-Frobenius theory that there exists an invariant class 
$\eta \in H^{1,1}_{psef}(X)$ associated to the spectral radius 
$\l=\r_{f^*} \geq 1$ of $f^*|_{\H11(X)}$.  When $\r_{f^*}>1$
(equivalent [GR] to saying that $f$ has positive entropy), it is
reasonable to hope that the associated current $T_\eta$ will itself be
positive.  In the next section, we pursue the construction of $T_\eta$,
and this hope in particular, for a much larger class of maps.

\section{Green's functions for meromorphic maps}

Let $X$ be a compact K\"ahler manifold of dimension $k$.  When 
$f:X \rightarrow X$ is merely meromorphic (i.e. rational, when $X$ is
projective), the construction of dynamical Green 
currents is a more delicate task, due to the presence of points of 
indeterminacy.  Nevertheless, Green currents have been constructed in some 
particular meromorphic cases (see [S] for the
case $X=\P^k$, [DF] for the case of birational surface maps, [G 1] for the case
of Hirzebruch surfaces, and [G 2] for a slightly more general context). 
In this section we use ideas of [BD] to provide a very general construction.

We let $I_f$ denote the indeterminacy locus, i.e. the set of points at which 
$f$ is not holomorphic. This is an analytic subset of $X$ of codimension 
$\geq 2$. We let $\Gamma_f \subset X \times X$ denote the graph of $f$
and $\tilde{\Gamma}_f$ denote a desingularization of it. We have a
commutative diagram
$$
\begin{array}{ccccc}
\text{ } & \text{ } & \tilde{\Gamma}_f & \text{ } & \text{ } \\
\text{ } & \stackrel{\pi_1}{\swarrow} & \text{ } &
\stackrel{\pi_2}{\searrow} & \text{ } \\
X & \text{ } & \stackrel{f}{\longrightarrow} & \text{ } & X
\end{array}
$$
where $\pi_1,\pi_2$ are holomorphic maps. We always assume that
$f$ is dominant, i.e. that its jacobian determinant does not vanish identically
(in any coordinate chart). 

Given a smooth closed real form $\o$ of bidegree $(1,1)$ on $X$, we set
$f^* \o:=(\pi_1)_* (\pi_2^* \o)$, where we push $\pi_2^* \o$ forward by 
$\pi_1$ as a current. Observe that $f^* \o$ is actually a form with 
$L^1_{loc}$-coefficients which coincides with the usual smooth pull-back 
$(f_{|X \setminus I_f})^* \o$ in $X \setminus I_f$.  Thus the definition 
does not depend on the choice of desingularization. Also $f^*$
preserves boundaries and thus induces an action
$f^*:\H11(X)\to\H11(X)$ given by
$$
\{\o\} \mapsto \{ f^* \o \}.
$$

Our formula for $f^*\omega$ may also be applied to pull back
(differences of) positive closed $(1,1)$ currents $S$: 
given $S\geq 0$, one uses the construction described at the end of the previous
section to define $\pi_2^*S$.  Then, as in the holomorphic case, 
$f^* S := \pi_{1*}\pi_2^* S$ is also a positive closed $(1,1)$
current.  It follows again that $f^*$ preserves the pseudoeffective
cone and that there exists $\eta\in H^{1,1}_{psef}(X)$ such that $f^*\eta
= \l\eta$ where $\l =\r_{f^*} \geq 1$ is the spectral radius of 
$f^*|_{\H11(X)}$.
  
An argument of M.Gromov [Gr] implies that $f$ has zero
topological entropy when $\r_{f^*}=1$. In the sequel we assume to the
contrary that $\r_{f^*}>1$.  Any smooth form $\omega$ representing $\eta$ may be
written as a difference of positive forms.  Hence we can iterate $f^*$
and try to construct a Green current $T_{\eta}=
\lim_{n\to\infty} f^{n*}\omega$ associated to $\eta$.

We immediately face some problems. First, the action $f^*|_{\H11(X)}$ is 
not necessarily compatible with the dynamics: it might happen, as in the
following example, that $(f^n)^*$ is different from $(f^*)^n$ for some 
$n \in \N$.

\begin{exa}
The polynomial endomorphism of $\C^2$
$$
(z_1,z_2) \mapsto (z_1 z_2,z_1).
$$
extends to a meromorphic endomorphism $f:\P^2\to\P^2$ of the complex
projective plane.  The extended map is given in homogeneous coordinates by
$$
f[z_0:z_1:z_2]=[z_0^2:z_1z_2:z_0 z_1],
$$
where $(z_0=0)$ denotes the line at infinity. Observe that the indeterminacy 
locus $I_f$ consists of the single point $[0:0:1]$. Since 
$\H11(\P^2) \simeq \R$ is one-dimensional, the linear action $f^*$ is
multiplication by $2$, so that $(f^*)^2$ is multiplication by $4$. On the
other hand, a simple computation shows that $(f^2)^*$ is multiplication by $3$.
\end{exa}

If one extends the polynomial map above to a meromorphic map $g$ on
$X=\P^1 \times \P^1$, one can check that $(g^*)^n=(g^n)^*$ for all $n \in \N$.
This motivates the following

\begin{defi}
The mapping $f:X \rightarrow X$ is said to be $1$-stable if
$(f^n)^*=(f^*)^n$ for all $n \in \N$.
\end{defi}

\begin{rqe}
The notion of 1-stable map has been studied by several authors in the past
decade, where it has been variously called \emph{generic} [FS 3],
\emph{minimally separating} [Di 1], \emph{algebraically stable} [S], or 
\emph{1-regular} [BK].

It was shown in [DF] that when $f$ is a birational surface map, one can always 
make a birational change of coordinates so that $f$ becomes 1-stable. It is an 
interesting, and probably quite difficult, open question to know whether this 
remains true for dominant 
2-dimensional meromorphic maps of `small' (i.e. less than 
$\rho(f^*|_{\H11(X)})$) topological degree.  Favre and Jonsson, in a recent
paper [FJ] concerning polynomial maps of $\C^2$, have proposed a very different
approach to issues concerning 1-stability. 
\end{rqe}

We assume from now on that $f:X \rightarrow X$ is a 1-stable meromorphic map.
Let $\eta \in H^{1,1}_{psef}(X)$ be such that $f^*\eta=\l \eta$,
with $\l=\r_{f^*}>1$. Let $\o$ be a smooth closed real $(1,1)$-form representing 
$\eta$.  By the $dd^c$-lemma again, there exists $\g_{\o} \in L^1(X)$ such that
\begin{equation}
\frac{1}{\l} f^* \o=\o+dd^c \g_{\o}.
\end{equation}
Since $f$ is 1-stable, we can pull this equation back by $f^n$ to get
$$
\frac{1}{\l^n} f^{n*} \o=\o+dd^c g_n^{\o}, \; \; 
\text{ where } \; \; g_n^{\o}:=\sum_{j=0}^{n-1} \frac{1}{\l^j} \g_{\o} \circ f^j.
$$

The second problem we face is that $\g_{\o}$ is not smooth when $f$ is 
meromorphic, so it is not obvious that the sequence $g_n^{\o}$ converges in 
$L^1(X)$.  This convergence is the content of our next result, which is a
slight refinement of Theorem A in the introduction.
Recall that a class $\eta \in H^{1,1}(X,\R)$ is {\it nef} if it is the limit 
of K\"ahler classes.

\begin{thm}
\label{construct}
If $\eta$ is nef, then the sequence $(g_n^{\o})$ converges in $L^1(X)$.
Let $g_{\o}$ be the limit, and define
$$
T_{\eta}:=\o+dd^c(g_{\o})=\lim_{n \rightarrow +\infty}
\frac{1}{\l^n} (f^n)^* \o.
$$
Then $T_\eta$ is a closed current satisfying $f^* T_{\eta}=\l T_{\eta}$.  If
$\lambda$ is a simple eigenvalue of $f^*$, then $T_\eta$ is positive.
\end{thm}

\begin{proof}
{\bf Step 1.}
The first step of the proof consists in showing that $\g_{\o}$ is 
bounded from above on $X$.  Rewriting $(1)$ in the 
desingularized graph yields
$$
dd^c \left( \g_{\o} \circ \pi_1 \right)=
\frac{1}{\l} \pi_1^* \pi_{1*}\pi_2^* \o - \pi_1^* \o.
$$
Since $\pi_1$ is a local isomorphism away from it's exceptional divisor
$\mathcal{E}(\pi_1)\subset \tilde{\Gamma}_f$, we have that 
$R:=\pi_1^* \pi_{1*}\pi_2^* \o-\pi_2^* \o$ is a closed current of
bidegree $(1,1)$ supported on ${\mathcal E}(\pi_1)$.  

\begin{lem} $R$ is positive.
\end{lem}

\proof
If $\omega$ (and hence $\pi_2^*\omega$) is a non-negative form, then so is 
$\pi_1^*\pi_1^*\pi_2^*\omega$.  Since $\pi_2^*\o$ agrees with the latter form
outside $\mathcal{E}(\pi_1)$ and does not charge $\mathcal{E}(\pi_1)$ itself,
it follows that $R$ is positive.  

Dropping the non-negativity assumption, we observe that $R$ depends only on
the cohomology class of $\omega$: if $\sigma = 
dd^c u$ is a cohomologically trivial $(1,1)$ form, then
$$
\pi_1^*\pi_{1*}\sigma - \sigma = dd^c (\pi_*\pi^* u - u) = dd^c 0 = 0.
$$
Since $\omega$ is nef, we may approximate $\omega$ \emph{in cohomology} by 
non-negative $(1,1)$ forms $\omega_j$.  Pullback and pushforward are continous
operators on currents, so
$R_j := \pi_1^*\pi_{1*} \pi_2^*\omega_j - \pi_2^* \omega_j$
converges weakly to $R$.  That is, $R$ is a limit of positive currents and
therefore positive.
\qed

Continuing with the proof of Theorem \ref{construct}, we have
$$
dd^c (\g_{\o} \circ \pi_1)=R+\pi_2^* \o-\pi_1^* \o,
$$
so that $\g_{\o} \circ \pi_1$ is a qpsh function on $\tilde{\Gamma}_f$.
Therefore $\g_{\o} \circ \pi_1$ is bounded from above on $\tilde{\Gamma}_f$,
as is $\g_{\o}$ on $X$.

{\bf Step 2.} 
It is now clear that we can add a constant to $\g_{\o}$ to
get $\g_{\o} \leq 0$. Therefore
$$
g_n^{\o}=\sum_{j=0}^{n-1} \frac{1}{\l^j} \g_{\o} \circ f^j
$$
is a decreasing sequence of $L^1$ functions.  We claim that it converges, i.e.
that $g_{\o}=\sum_{j \geq 0} \l^{-j} \g_{\o} \circ f^j$ belongs to $L^1(X)$. 

The claim can be established using tricky integration by parts as in [G 2] 
(pp 2377/78). Instead, we follow here a trick of N. Sibony (who treats the 
case $X=\P^k$ in [S]).  Since $\omega$ is nef, we can choose K\"ahler forms
$\omega_j$ whose classes converge to that of $\omega$.  The mass of $\omega_j$
is controlled by its cohomology class, so the sequence $(\omega_j)$
must accumulate on some positive closed current $\tilde \omega$ cohomologous 
to $\omega$.  

Considering Cesaro means of $\l^{-j}(f^j)^* \tilde \omega$ and extracting 
a limit, we can produce a positive closed $(1,1)$-current $\sigma$ which is 
cohomologous to $\o$ and satisfies $f^* \sigma=\l \sigma$.  Now 
$\sigma=\o+dd^c v$ for some qpsh function $v$ on $X$.  
Invariance of $\sigma$ allows us to arrange 
$$
\g_{\o}=v-\frac{1}{\l} v \circ f
$$
by adding a constant to $v$.  Pulling this equation back by $f^n$ then gives
$$
g_n^{\o}=v-\frac{1}{\l^n}v \circ f^n \geq v- \frac{\sup_X v}{\l^n}.
$$
Thus, $v \leq g_{\o} \leq 0$ and in particular $g_{\o} \in L^1(X)$.
The current $T_{\eta}:=\o+dd^c(g_{\o})$ clearly satisfies 
$f^* T_{\eta}=\l T_{\eta}$.  

{\bf Step 3.}  It remains to prove that when $\lambda$ is a simple eigenvalue 
of $f^*$, then $T_\eta$ is positive.  Taking 
$p=\dim_\C X$, we observe that it suffices to show that
$$
\pair{T_\eta}{\chi\sigma} \geq  0,
$$
where $\chi$ is a smooth non-negative cutoff function supported on a
coordinate chart $U\subset X$ and $\sigma$ is a positive $(p-1,p-1)$ form
that is constant with respect to coordinates on $U$.

\begin{lem}
If $\lambda$ is a simple eigenvalue of $f^*$, then
some subsequence of $(f^n_*(\chi\sigma)/\lambda^n)$ converges weakly to a 
positive closed $(p-1,p-1)$ current $S$.
\end{lem}

\begin{proof}
Note first that if $\omega_0$ is a K\"ahler form on $X$, then  
$0\leq \chi\sigma \leq C\omega_0^{p-1}$ for $C>0$ large enough.  Hence
assuming $\lambda$ is a simple eigenvalue of $f^*$, we obtain uniform control
on the mass of $\lambda^{-j} f^n_*(\chi\sigma)$ as follows.
$$
\int \omega_0\wedge f^n_*(\chi\sigma) \leq 
C \int \omega_0\wedge f^n_*\omega_0^{p-1}
= 
C\int f^{n*}\omega_0\wedge \omega_0^{p-1} \leq C'\lambda^n.
$$
In particular, the sequence $(f^n_*\,\chi\sigma/\lambda^n)$ has weak limit 
points.  These will be positive by continuity, so we need only 
show that they are also closed.  For this we employ an a well-known argument
of (see [BS]) to show that the mass of 
$\partial f^n_*(\chi\sigma)$ is no larger than $C\lambda^{n/2}$.  Specifically,
we let $\varphi$ be any real test $1$-form and estimate
\begin{eqnarray*}
|\pair{\partial f^n_* (\chi\sigma)}{\varphi}|
& = &
\left|\int f^{n*}\varphi\wedge d\chi\wedge \sigma \right| \\
& \leq &
\left(\int f^{n*}\varphi\wedge Jf^{n*}\varphi \wedge \sigma \right)^{1/2}
\left(\int d\chi\wedge d^c\chi \wedge \sigma \right)^{1/2} \\
& \leq &
C \left(\int f^{n*}\omega_0\wedge\sigma\right)^{1/2}
\left(\int d\chi\wedge d^c\chi \wedge \sigma \right)^{1/2}
\leq C'\lambda^{n/2}.
\end{eqnarray*}
Note that $J$ here is the complex structure operator on real cotangent
vectors. Moreover, all integrals may be interpreted as taking place away from 
the set $I(f)$ where $f^{n*}\varphi$ might be singular.  Having established
the desired control on $\partial f^n_*(\chi\sigma)$, we are done. 
\end{proof}

Using the subsequence from the lemma, we have
$$
\pair{T_\eta}{\chi\sigma} = \lim_{j\to\infty} \lambda^{-n_j} 
\pair{f^{n_j*}\omega}{\chi\sigma}
=
\lim_{j\to\infty} \lambda^{-n_j} \pair{\omega}{f^{n_j}_*(\chi\sigma)} 
= 
\pair{\omega}{S}
= \eta \cdot\{S\} \geq 0,
$$
where the last inequality comes from the assumption that $\eta$ is nef.
\end{proof}
\vskip.1cm

It remains to understand when our hypotheses are satisfied.
When $\dim_{\C} X=2$, the cone $H^{1,1}_{nef}(X,\R)$ is preserved by $f^*$,
so the invariant class $\eta$ is automatically nef (Proposition 1.11 in
[DF]).  Moreover, $\lambda^2$ is never less than the topological degree of
$f$; when it is strictly larger, it is a simple eigenvalue of $f^*$ (Remark
5.2 in [DF]).

\begin{cor}
\label{dim2cor}
If $\dim_{\C} X=2$, the sequence $(g_n^{\o})$ always converges in $L^1(X)$.
If $\lambda^2$ exceeds the topological degree of $f$, then the associated
closed invariant current $T_\eta := \omega + dd^c g_\omega$ is positive.
\end{cor}

For the rest of the paper, we focus exclusively on the case $\dim_\C X = 2$
of maps on complex surfaces.
Recall from [S] that a point $x \in X$ is said to be {\it normal} if there
exist neighborhoods $U$ of $x$ and $V$ of $I_f$, such that $f^nU \cap
V=\emptyset$ for all $n \in \N$.  The set of normal points is denoted by 
${\mathcal N}_f$: this is the set of points which remain `locally uniformly' 
away from the indeterminacy locus under iteration.  The proof of Theorem 1.1
applies straightforwardly to show that $g_{\o}$ is H\"older continuous
in ${\mathcal N}_f$ (this is Theorem 1.7.1 in [S]). 
Complex H\'enon mappings are polynomial automorphisms of $\C^2$ which extend to
$\P^2$ as 1-stable maps of positive entropy. For such mappings the set of 
normal points is ${\mathcal N}_f=\P^2 \setminus I_f$, hence the dynamical Green
function is H\"older continuous off the indeterminacy locus. This result was first proved in [FS 1].

\section{Sub-H\"older continuity}

The set ${\mathcal N}_f$ might well be empty for a given meromorphic map $f$,
and $g_{\o}$ can be very discontinuous in general (see e.g. example 1.11 in 
[GS] and example 6.1 below).  In this section we consider some families of 
rational surface mappings
that permit weaker, though still `H\"older-like' control on the modulus of 
continuity of the dynamical Green's function.
$$
g_{\o}:=\sum_{j=0}^{+\infty} \frac{1}{\l^j} \g_{\o} \circ f^j.
$$
Of course this can be done only off the extended indeterminacy locus,
$$
{\mathcal I}_f:=\overline{ \cup_{n \geq 0} I_{f^n} },
$$
since $g_{\o}$ usually has positive Lelong number at every point of $I_{f^n}$.

H\"older continuity of $g_{\o}$ at $p$ requires that the orbit of $p$
uniformly avoid the indeterminacy locus (normal points).  Weaker kinds of 
continuity of $g_{\o}$ can be established by simply requiring that $f^n(p)$ 
not approach $I_f$ too rapidly: see [FG], [G 1], [GS] for the case of 
weakly-regular polynomial endomorphisms of $\C^k$; and [Di 2] for birational 
maps of $\P^2$ that are {\it separating}, i.e. such that
${\mathcal I}_f \cap {\mathcal I}_{f^{-1}}=\emptyset$.  

Here we present a unified approach to estimating the modulus of continuity 
of $g_{\o}$ in $X \setminus {\mathcal I}_f$.  It applies to a class of 
rational maps large enough to encompass both weakly-regular endomorphisms
and separating birational surface maps.  Our main dynamical assumption 
is as follows. 
\emph{There exists $C>1$ and $\b \in [1,\l[$ such that}
\begin{equation}
\frac{1}{C} \left[ d(x,\I_f) \right]^{\b} \leq d(fx, \I_f) , \; \; 
\text{ for all } x \in X \setminus {\mathcal I}_f.
\end{equation}
Note that this estimate is stated in terms of the extended indeterminacy set
$\I_f$ rather than just the indeterminacy set $I_f$.  This is because $I_f$
is rarely invariant under $f^{-1}$, whereas one always has 
$f^{-1}(\I_f - I_{f^{-1}})\subset \I_f$.

We shall rely on three further estimates, all of which hold independent of the
above assumption.  The first gives us pointwise
control on $\g_{\o}$:
\begin{equation}
\label{logbound}
\g_{\o}(x) \geq C \log d(x,I_f) + C'
\end{equation}
This follows from Proposition 1.2 in [BD] which, despite the birational
context of that paper, remains valid for arbitrary meromorphic surface
maps.  The other two estimates are local bounds on the Lipschitz constants
of $f$ and $\g_{\o}$.  Namely, one can check
by computing in local charts that there exist $m_1,m_2>0$ such that for all
$x,y \in X \setminus I_f$,
\begin{equation}
|\g_{\o}(x)-\g_{\o}(y)| \leq \frac{C d(x,y)}{[d_{I_f}(x,y)]^{m_1}}
\; \text{ and } \; 
d(f(x),f(y)) \leq \frac{C d(x,y)}{[d_{I_f}(x,y)]^{m_2}},
\end{equation}
where $d_{I_f}(x,y):=\min \{d(x,I_f),d(y,I_f) \}$ is the distance from
the pair $\{x,y\}$ to the indeterminacy locus.  We similarly denote distance
to $\I_f$ by $d_{\I_f}(x,y)$.  For convenience, we take
the constant $C>0$ to be the same in (2),(3), and (4).
It follows directly from (4) that
\begin{equation}
|\g_{\o} \circ f^j(x) - \g_{\o} \circ f^j(y))| \leq
\frac{C^{j+1} d(x,y)}{ [d_{I_f}(f^jx,f^jy)]^{m_1} 
\Pi_{l=0}^{j-1} [d_{I_f}(f^lx,f^ly)]^{m_2}}.
\end{equation}
It follows from (2) and (4) and the fact that $d_{\I_f}\leq d_{I_f}$ that
\begin{equation}
\frac{1}{d_{\I_f}(f^jx,f^jy)} \leq \frac{C^{1+\b+\cdots+\b^{j-1}}}{[d_{\I_f}(x,y)]^{\b^j}}.
\end{equation}

We will use these bounds to obtain control on the modulus of continuity of 
$g_{\o}$.  We treat the cases $\b=1$ and $\b>1$ separately since they are 
quite different.

\subsection{The case $\b=1$}
Our aim here is to prove the following

\begin{thm}
Let $f:X \rightarrow X$ be a 1-stable map which satisfies (2) with $\b=1$.
Then there exists $\a>0$ such that for all $x,y \in  X \setminus {\mathcal I}_f$,
$$
|g_{\o}(x)-g_{\o}(y)| \leq C_{x,y} \exp (-\a \sqrt{|\ln d(x,y)|} ),
$$
where $(x,y) \mapsto C_{x,y}>0$ is locally uniformly bounded in $X \setminus {\mathcal I}_f$.
\end{thm}

We need the following elementary lemma whose proof is left to the reader.

\begin{lem}
Fix $\a \in ]0,1[$ and set, for $ 0 \leq t \leq 1$,
$$
h_{\a}(t):=\exp(-\a \sqrt{|\ln t|}).
$$
Then for all $t \in [0,1]$ and for all $A \geq 1$,
$$
0 \leq t \leq e \,  h_{\a}(t)
\text{ and }
0 \leq h_{\a}(A t) \leq \exp(\a \sqrt{\ln A}) h_{\a}(t).
$$ 
\end{lem}

\begin{proof.3.1} 
Let $x,y\in X-\I_f$ be given.
Since $\b=1$, we infer from (3) and (6) that
$|\gamma\circ f^j(x)| \leq Cj$
where $C$ depends only on $d(x,\I_f)$.  In particular, if $t>1$ then
$$
\frac{|\gamma\circ f^j(x)-\gamma\circ f^j(y)|}{t^j} \leq M
$$
for all $j\geq 0$ and some constant $M$ depending on $t$ and $d_{\I_f}(x,y))$.
Moreover, from (5) and (6), we obtain the alternative upper bound
$$
\frac{|\g_{\o} \circ f^j(x) - \g_{\o} \circ f^j(y)|}{t^j} \leq
\frac{C^{m_2j^2+m_1j+j+1}}{t^j[d_{\I_f}(x,y)]^{m_1+m_2j}} d(x,y).
$$
It follows therefore from lemma 3.2 that 
$$
\frac{|\g_{\o} \circ f^j(x) - \g_{\o} \circ f^j(y)|}{Mt^j} \leq 
e \exp \left(\a \sqrt{ (m_2j^2\ln C +Aj + B)} \right)
h_{\a}(d(x,y)),
$$
where $A$ and $B$ depend on $t$ and $d_{\I_f}(x,y)$.  Thus
$$
|g_{\o}(x)-g_{\o}(y)| \leq \sum_{j \geq 0} \frac{Mt^j}{\l^j} 
\frac{|\g_{\o} \circ f^j(x)-\g_{\o} \circ f^j(y)|}{Mt^j} 
\leq C_{x,y,t} h_{\a} \circ d(x,y),
$$
where the series defining $C_{x,y,t}$ converges as soon as 
$\a<\ln(\l/t) [m_2 \ln C]^{-1/2}$ since it is comparable to
$$
\sum_{j \geq 0} \frac{t^j}{\l^j} \exp \left[ j \a \sqrt{m_2 \ln C} \right] 
<+\infty.
$$
Observe also that the dependence of $C_{x,y}$ on $(x,y)$ only involves 
$\ln d_{\I_f}(x,y)$, hence it is bounded on compact subsets of 
$X \setminus {\mathcal I}_f$.
\end{proof.3.1}
\vskip.2cm

We now want to provide examples of rational mappings satisfying the
assumptions of Theorem 3.1.  Observe first that any {\it separating} 
birational 1-stable self-map of a compact K\"ahler surface $X$ satisfies (2): 
this was observed by the first author who proved Theorem 3.1 in this context 
(see Theorem 5.3 in [Di 2]). The reader will find several examples of such 
birational mappings in [Di 1]. Note that the proof given above greatly 
simplifies the proof given in [Di 2]. It also applies to non birational 
mappings, as the following example shows.

\begin{exa}
Consider the meromorphic map $f:\P^2\to\P^2$, given in local
coordinates on $\C^2$ by
$$
f:(z_1,z_2) \in \C^2 \mapsto (P(z_1),Q(z_1)+R(z_2)) \in \C^2,
$$
where $P,Q,R$ are polynomials of degree $p,q,\l$ with $\l=pq >\max(p,q)$.
The indeterminacy locus consists of the single point
$$
I_f=[0:0:1]=\{z_0=z_2=0\},
$$
where $\{z_0=0\}$ denotes the line at infinity, written in homogenous coordinates.
Observe that $f$ contracts the line at infinity to the superattracting fixed point
$[0:1:0] \notin I_f$. Therefore $f$ is 1-stable on $\P^2$ and
${\mathcal I}_f=I_f$. The map $f$ is an example of weakly-regular polynomial endomorphism
of $\C^2$ (see [GS]). It corresponds to the critical case of mappings considered in
[G 1], since the topological degree of $f$ 
$$
d_t(f)=pq=\l
$$
coincides with its first dynamical degree $\l>1$.  Computations similar to
those in [G 1] show that $f$ satisfies (2) with $\b=1$.
\end{exa}

\subsection{The case $\b>1$}

\begin{thm}
Assume $f$ satisfies (2) with $1<\b<\l$.  Fix 
$\a \in [0,\frac{\ln \l}{\ln \b}-1[$.
Then for all $x,y \in X \setminus {\mathcal I}_f$,
$$
|g_{\o}(x)-g_{\o}(y)| \leq \frac{C_{\a}(x,y)}{1+|\ln d(x,y)|^{\a}},
$$
where $(x,y) \mapsto C_{\a}(x,y)$ is locally uniformly bounded in 
$X \setminus {\mathcal I}_f$.
\end{thm}

We need the following elementary lemma whose proof is left to the reader.

\begin{lem}
Set $\f_{\a}(t):=[1+|\ln t|^{\a}]^{-1}$, for $ 0 \leq t \leq 1$.

Then there exists $C'>1$ independent of $\a$ such that for all $t \in [0,1]$
and for all $A \geq e$,
$$
0 \leq t \leq C' \f_{\a}(t) \; \;
\text{ and } \; \;
0 \leq \f_{\a}(A t) \leq C' (\ln A)^{\a} \f_{\a}(t).
$$
\end{lem}

\begin{proof.3.4}
Let $x,y\in X-\I_f$ be given.
Since $\b>1$ it follows from (3) and (6) that
$$
|\g_\o\circ f^j(x) - \g_\o\circ f^j(y)| \leq \max\{|\g_\o \circ f^j(x)|,|\g_\o \circ f^j(y)|\}
\leq M\beta^j,
$$
for all $j\in\N$ and some constant $M$ depending on $d_{\I_f}(x,y)$.  From (5) 
and (6) we have the alternative bound
$$
|\g_{\o} \circ f^j(x) - \g_{\o} \circ f^j(y)| \leq 
\left[ \frac{C_2}{d_{\I_f}(x,y)^{m_3}} \right]^{\b^j} d(x,y),
$$
where $m_3=m_1+m_2/(\b-1)$ and $\ln C_2=[1+m_1/(\b-1)+m_2/(\b-1)^2] \ln C$.
We infer using lemma 3.5,
$$
\frac{|\g_{\o} \circ f^j(x) - \g_{\o} \circ f^j(y)|}{M\beta^j} \leq 
C\b^{\a j} \f_{\a} \circ d(x,y)
$$
for all $j\in\N$.  Once again $C$ depends on $x$ and $y$ only through 
$d_{\I_f}(x,y)$.  The factor of $M\beta^j$ ensures that the right side remains 
bounded above by $1$; i.e. it guarantees that the hypothesis of lemma 3.5 is 
satisfied.

Now we conclude
$$
|g_\o(x)-g_\o(y)| \leq \sum_{j=0}^\infty M\left(\frac\beta\lambda\right)^j
                  \frac{|\g_\o\circ f^j(x)-\g_\o \circ f^j(y)|}{M\beta^j}
                  \leq C \f_\alpha(d(x,y))\sum_{j=0}^\infty 
                  \left(\frac{\beta}{\lambda}\right)^j\beta^{\alpha j}. 
$$
The sum on the right side converges as soon as 
$\beta^{\a+1}< \lambda$, i.e. as soon as 
$\alpha < \ln\lambda/\ln\beta -1$.
\end{proof.3.4}

\vskip.1cm
Weakly-regular polynomial endomorphisms of $\C^k$ as considered in [FG], [G 1], [GS]
provide several examples of mappings which satisfy the assumptions of Theorem 3.4.
Let us recall that such control on the modulus of continuity of $g_{\o}$
yields integrability of $\log^+ ||Df||$ with respect to special invariant measures,
as well as estimates on the pointwise dimension of these measures 
(see [S], [Di 2]).

\section{Sobolev regularity}

It can happen that the extended indeterminacy locus ${\mathcal I}_f$ is very
large (see section 6 for an example where ${\mathcal I}_f=X$).  Our aim in 
this section is to consider weaker regularity properties for the 
functions $g_{\o}$ which can hold even across points of indeterminacy.
To simplify {\bf we restrict ourselves to the case where
$f$ is a bimeromorphic surface map}, i.e. $\dim_{\C} X=2$ and there
is a meromorphic map $f^{-1}:X \rightarrow X$ such that 
$f \circ f^{-1}=\text{Id}$.

In this case it was proved by the first author and C.Favre [DF] that one can
always make a birational change of coordinates so that $f:X \rightarrow X$ 
becomes $1$-stable. Moreover the spectral radius $\l$ of 
$f^*:H^{1,1}_{\R}(X) \rightarrow H^{1,1}_{\R}(X)$ is a simple eigenvalue as 
soon as $\l>1$.  Since $f^*$ preserves the pseudoeffective cone and $f_*$ is
intersection adjoint to $f^*$, there are classes 
$\{ \o^{\pm} \} \in H^{1,1}_{psef}(X)$, unique up to positive multiple, 
such that
$$
f^* \{\o^+ \}=\l \{ \o^+ \} \text{ and }
(f^{-1})^* \{\o^- \}=\l \{\o^- \},
$$
By Theorem 2.4 applied to $f$ and $f^{-1}$, there are positive closed currents
$$
T^{\pm}=\o^{\pm}+dd^c g^{\pm},
\text{ where } g^{\pm}=\sum_{j \geq 0} \frac{1}{\l^j} \g^{\pm} \circ f^{\pm j}
$$
and the functions $\g^{\pm} \in L^1(X)$ satisfy 
$\l^{-1}(f^{\pm})^* \o^{\pm}=\o^{\pm}+dd^c \g^{\pm}$.
We shall assume moreover that
\begin{equation}
\label{minimality}
\{ \o^+ \} \cdot f(p) >0 \; \text{ and } \;
\{ \o^- \} \cdot f^{-1}(q)>0,
\end{equation}
for all $p\in I_f$, $q\in I_{f^{-1}}$. This can always be arranged (see 
Proposition 4.1 in [BD]).
\vskip.1cm

The main result of this section identifies geometric conditions equivalent 
to the statement that the gradients $\nabla g^{\pm}$ belong to $L^2(X)$.
Since qpsh functions are always in $L^1(X)$, this
is equivalent to saying that $g^+$ and $g^-$ belong to the Sobolev space
$W^{1,2}(X)$.

\begin{thm}
Let $f:X \rightarrow X$ be a birational map with $\l=\r_{f^*}>1$. Then the following conditions
are equivalent:
\begin{enumerate}
\item $\nabla g^+ \in L^2(\o^- \wedge \Omega)$;
\vskip.2cm
\item for all $p \in I_f$, $\sum_{n \geq 0}  
\frac{\g^- \circ f^{-n}(p)}{\l^{2n}} >-\infty$;
\vskip.2cm
\item $\sum_{n \geq 0} \l^{-2n} \log {\rm dist}(I_{f^{-1}},f^{-n}I_f) >-\infty$.
\end{enumerate}
\end{thm}

Observe that when $\o^-$ is K\"ahler, then (1) means precisely that $g^+$ has 
gradient in $L^2(X)$ with respect to the Lebesgue measure. The class $\{ \o^-
\}$ is automatically K\"ahler when 
for instance $X=\P^2$. This condition should be compared to the slightly stronger
condition studied in [BD]:  $g^-$ is finite at
each point of $I_f$. This is equivalent to
$$
\sum_{n \geq 0} \frac{\g^- \circ f^{-n}(p)}{\l^{n}} >-\infty
$$
for all $p \in I_f$. As in [BD], these equivalent conditions are symmetric in
$f$ and $f^{-1}$; i.e. one can interchange the roles played by $f$ and its
inverse $f^{-1}$ and obtain three further equivalent conditions.

\begin{proof}
The equivalence between (2) and (3) follows from the fact that under the 
assumption \eqref{minimality},
$\g^-$ is a smooth function in $X \setminus I_{f^{-1}}$ with logarithmic 
singularities
at points of indeterminacy of $f^{-1}$. More precisely, there exists constants
$A,B,A',B'>0$ such that
$$
A \log \text{dist}(x,I_{f^{-1}}) -B \leq \g^-(x) \leq 
A' \log \text{dist}(x,I_{f^{-1}}) -B'.
$$
We refer the reader to [BD] for a proof of this fact.

It is a simple exercise to check that condition (2) is equivalent to
the finiteness of the sum
$\sum_{n \geq 0} \l^{-n} g_n^-(p)$
at all points $p \in I_f$.
Therefore the equivalence between (1) and (2) is a consequence of the next lemma.
\end{proof}

\begin{lem}
$$
\int_X dg_n^+ \wedge d^c g_n^+ \wedge \o^-=\sum_{p \in I_f} \sum_{j = 0} ^{n-1} c_j(p)
\frac{|g_j^-(p)|}{\l^j}+O(1),
$$
where the constants $c_j(p)$ are positive and uniformly bounded away from $0$
and $\infty$.
\end{lem}

\begin{proof}
Set $T_n^+:=\l^{-n}(f^n)^* \o^+$. For $0 \leq j \leq n-1$ we have
\begin{equation}
\int (\g^+ \circ f^j)\, T_n^+ \wedge \o^-=\int \g^+ T_{n-j}^+ \wedge 
\frac{f^j_*\o^-}{\l^j}
=\int \g^+ T_{n-j}^+ \wedge \o^-+A_{n,j},
\end{equation}
where
$$
A_{n,j}:=\int \g^+ T_{n-j}^+ \wedge dd^c g_j^-=
\int \left[ \frac{1}{\l}f^* \o^+ -\o^+ \right] \wedge g_j^- T_{n-j}^+.
$$
Observe that the currents $f^* \o^+, T_{n-j}^+$ both have positive Lelong 
numbers at points in $I_f$, thus
$$
\frac{1}{\l}f^* \o^+ \wedge T_{n-j}^+=\sum_{p \in I_f} c_{n-j}(p) \d_p
+\frac{1}{\l^2}f^* \left( \o^+ \wedge T_{n-j-1}^+ \right),
$$
where $\d_p$ denotes the Dirac mass at point $p$ and the $c_j$'s are positive 
constants.  Observe that $f^* \o^+ \wedge T_j \rightarrow f^*\o^+ \wedge T^+$,
since $g_j^+$ decreases towards $g^+$.  Since $T^+$ has positive Lelong number 
at all points of indeterminacy, it follows that the measure 
$f^* \o^+ \wedge T^+$ has a Dirac mass $c_{\infty}(p)>0$ at each point 
$p \in I_f$.  Therefore $c_j(p) \rightarrow c_{\infty}(p)>0$.  In particular, 
the sequences $(c_j(p))$ are uniformly bounded away from zero and infinity. We 
infer
$$
A_{n,j}=\sum_{p \in I_f} c_{n-j}(p) g_j^-(p)+\frac{1}{\l^2} 
\int (g_j^- \circ f^{-1})\, \o^+ \wedge T_{n-(j+1)}^+
-\int g_{j}^- \o^+ \wedge T_{n-j}^+.
$$
Observe that $\l^{-1} g_j^- \circ f^{-1}=g_{j+1}^- -\g^-$. Thus multiplying (8) by
$\l^{-j}$ and summing from $j=0$ to $n-1$ yields
$$
\int g_n^+ T_n^+ \wedge \o^-=\sum_{p \in I_f} 
\sum_{j=0}^{n-1} c_{n-j}(p) \frac{g_j^-(p)}{\l^j} +
\int \frac{g_n^-}{\l^n} \o^+ \wedge \o^+ -\int \g^- T^+ \wedge \o^+ +M_n,
$$
where
$$
M_n:=\sum_{j=0}^{n-1} \frac{1}{\l^j} \int  \g^+ T_{n-j}^+ \wedge \o^+
-\sum_{j=0}^{n-1} \frac{1}{\l^{j+2}} \int \g^- \o^+ \wedge T_{n_j-1}^+
$$
is a bounded sequence. Note to conclude that
\begin{eqnarray*}
\int (-g_n^+) T_n^+ \wedge \o^-&=& \int (-g_n)^+ dd^c g_n^+ \wedge  \o^-+\int (-g_n)^+ \o^+ \wedge \o^+ \\
&=& \int d g_n^+ \wedge d^c g_n^+\wedge\o^- + O(1).
\end{eqnarray*}
\end{proof}

The function $g^+$ (resp. $g^-$) typically has positive Lelong number at every 
point of $I_f^{\infty}:=\cup_{n \geq 0} I_{f^n} = \cup_{n\geq 0} f^{-n}(I_f)$ 
(resp. $I_{f^{-1}}^{\infty}:=\cup_{n \geq 0} I_{f^{-n}} = 
\cup_{n\geq 0} f^n(I_{f^{-1}})$). The closer $f^{-n} I_{f}$ is to
$I_{f^{-1}}$, the stronger the singularity of $g^-$ at these points.
The symmetric condition $\nabla g^{\pm} \in L^2(X)$ is thus a quantitative way 
of saying that the sets $I_f^{\infty}$ and $I_{f^{-1}}^{\infty}$ stay away 
from each other (recall that the condition of 1-stability means precisely that
$I_f^{\infty} \cap I_{f^{-1}}^{\infty}=\emptyset$).  We will see in section 6 
that this condition is often but not always satisfied.

\section{A canonical invariant measure}

We consider here, as in the previous section, a compact K\"ahler surface $X$
equipped with a K\"ahler form $\Omega$, and a 1-stable bimeromorphic map
$f:X \rightarrow X$ such that $\l=\rho_f>1$.
As before, we let $T^{\pm}=\o^{\pm}+dd^c g^{\pm}$ denote the 
positive closed $(1,1)$-currents invariant under $f^{\pm 1}$.
We normalize so that
$$
\{T^+ \} \cdot \{T^- \}=\{ \Omega \} \cdot \{T^+ \}=
\{\Omega \} \cdot \{T^- \}=1.
$$
Our aim here is to define and study the measure
$$
\mu_f:=\text{``} T^+ \wedge T^- \text{''}.
$$

\subsection{Definition of the canonical measure $\mu_f$}
It is well known that one cannot always define the wedge product of two 
positive closed currents. When $g^+$ is integrable with respect to the trace
measure of $T^-$, i.e. $g^+ \in L^1(T^- \wedge \Omega)$, the current $g^+ T^-$ 
is well defined, and we can set
$$
\mu_f:=\o^+ \wedge T^-+dd^c (g^+ T^-).
$$
Observe that the condition is symmetric, as follows from Stokes theorem:
$$
\int_X (-g^+) T^- \wedge \Omega=\int_X (-g^-) T^+ \wedge \Omega +
\int_X (g^- \o^+-g^+ \o^-) \wedge \Omega.
$$

When the potentials $g^{\pm}$ have gradients in $L^2(X)$, it follows from 
the Cauchy-Schwarz inequality that $g^+ \in L^1(T^- \wedge \Omega)$, since
\begin{eqnarray*}
\lefteqn{
0 \leq \int -g^+ T^- \wedge \Omega=\int -g^+ \o^- \wedge \Omega
+\int dg^+ \wedge d^cg^- \wedge \Omega} \\
&\leq& \! \! \! \! \int -g^+ \o^- \wedge \Omega+
\left( \int dg^+ \wedge d^c g^+ \wedge \Omega \right)^{1/2}
\left( \int dg^- \wedge d^c g^- \wedge \Omega \right)^{1/2}
\! \! \! \! <+\infty.
\end{eqnarray*}

It may happen, however, that $g^+ \in L^1(T^- \wedge \Omega)$ while
$\nabla g^+ \notin L^2(X)$ (see example 6.2).
We know of no example for which the function 
$g^+$ is not integrable with respect to the trace measure of $T^-$.  Hence
the following:

\begin{ques}
Is the condition $g^+ \in L^1(T^- \wedge \Omega)$ always satisfied ?
\end{ques}

We now derive a criterion which will allow us to check the condition
$g^+ \in L^1(T^- \wedge \Omega)$ for some birational mappings.

\begin{pro}
Let $S$ be a positive closed (1,1)-current on $X$, whose cohomology class
$\{S\}$ is K\"ahler. Assume
\begin{enumerate}
\item $g^- \in L^1(S \wedge \Omega)$, so that the measure $S \wedge T^-$ is well defined; 
\item $g^+ \in L^1(S \wedge T^-)$.
\end{enumerate}
\hskip.2cm  Then $g^+ \in L^1(T^- \wedge \Omega)$.
\end{pro}

\begin{proof}
Let $\theta_S \geq \e \Omega$ be a K\"ahler form cohomologous to $S$.
Let $\f_S \in L^1(X)$ be a qpsh function such that
$S=\theta_S+dd^c \f_S$. We can assume without loss of generality that
$\f_S \leq 0$ and $\e=1$. Then
$$
0 \leq \int_X -g^+ T^- \wedge \Omega \leq \int_X -g^+ T^- \wedge \theta_S
= \int_X -g^+ T^- \wedge S+\int_X g^+ T^- \wedge dd^c \f_S.
$$
The next to last integral is finite by assumption. The last one is finite by 
Stokes theorem:
\begin{eqnarray*}
\int g^+ T^- \wedge dd^c \f_S 
& = & \lim_{n\to\infty} \int g_n^+ T^- \wedge dd^c \f_S 
= \lim_{n\to\infty} \int -\f_S T^-\wedge (-dd^c g_n^+) \\
& \leq & \int_X (-\f_S) T^- \wedge\o^+
< \infty.
\end{eqnarray*}
The first inequality holds because 
$-dd^c g_n^+ = \omega^+ - T^+ \leq \omega^+$.  The second holds
because $\o^+ \leq C \Omega$ and 
$\f_S \in L^1(T^- \wedge \Omega) \Leftrightarrow g^- \in L^1(S \wedge \Omega)$.
\end{proof}

We will use this criterion in section 6.2 when $S=[V]$ is the current of integration 
along an invariant irreducible curve $V$.

\subsection{Dynamical properties of $\mu_f$}

We assume in the sequel that $g^+ \in L^1(T^- \wedge \Omega)$ so that
$\mu_f=T^+ \wedge T^-$ is well defined.

\begin{thm}
The measure $\mu_f$ does not charge the indeterminacy locus $I_f$. It is
an invariant probability measure. 
\end{thm}

\begin{proof}
The current $\mu_f$ is a positive measure. This can be seen by locally regularizing
the qpsh function $g^+$. It is a
probability measure by our choice of normalization,
and it is the weak limit of the measures
$$
\mu_n:=\frac{1}{\l^n} (f^n)^* \o^+ \wedge T^-.
$$
Let $\chi$ be a test function. Observe that $T^-$ does not charge curves
(it has zero Lelong number at each point in $X \setminus I_{f^{-1}}^{\infty}$,
see [DF]), hence neither does $\mu_n$. It follows therefore from the change
of variables formula and the invariance $f_*T^-=\l T^-$, that
$$
\langle \mu_n, \chi \rangle=\frac{1}{\l^{n+1}} \langle (f^n)^* \o^+, \chi f_*T^- \rangle
=\langle \mu_{n+1}, \chi \circ f \rangle.
$$
We infer, if $\mu_f(I_f)=0$, that $\langle \mu_f,\chi\rangle=\langle \mu_f, \chi \circ f \rangle$,
hence $\mu_f$ is invariant.

It remains to prove that $\mu_f$ does not charge any point $p\in I_f$.  
Since $f$ is $1$-stable, we may assume that 
$p \notin I_{f^{-1}}^{\infty}$.  If $p = f^{-N}(p)$ is periodic,
then $g^-$ is finite at $p$.  Hence it follows from Proposition 5.4 below 
that $\mu_f(\{p\})=0$.

If on the other hand, $p$ is not periodic, we can fix $N >>1$ and
choose $r=r_N>0$ such that $f^{-j}B(p,r) \cap f^{-k}B(p,r)=\emptyset$ for
$0 \leq j,k \leq N$, $j \neq k$. Let $\chi$ be a test function such that
$0 \leq \chi \leq 1$, $\chi \equiv 1$ near $p$, and
$\supp \, \chi \subset B(p,r)$. Then
$$
0 \leq \p:=\sum_{j=0}^{n-1} \chi \circ f^j \leq 1,
$$
since the functions $\chi \circ f^j$ have disjoint supports.
Set 
$$
R_n:= (\p \circ f^n) T^+=\sum_{j=0}^{N-1} (\chi \circ f^{j+n}) T^+.
$$
It follows from the extremality of $T^+$ that
the positive currents $(\chi \circ f^{j+n}) T^+$ converge
to $c_{\chi} T^+$, where $c_\chi=\int \chi d\mu_f$ (see [BD], [G 1]).
Thus $R_n \rightarrow N c_{\chi} T^+$. Now
$0 \leq R_n \leq T^+$ since $0 \leq \p \leq 1$, hence
$$
0 \leq \mu_f(\{p\}) \leq c_{\chi}=\int \chi d\mu_f \leq \frac{1}{N}.
$$
Since $N$ was arbitrary, we conclude that 
$\mu_f(\{p\})=0$.  Thus $\mu_f$ does not charge the indeterminacy locus.
\end{proof}

\begin{pro}
\label{finiteg}
Fix $p \in I_f$. Then $g^-(p)$ is finite if and only if 
$\log {\rm dist}(\cdot,p) \in L^1(\mu_f)$.
In particular $g^-$ is finite on $I_f$ (the [BD] condition described after
Theorem 4.1) if and only if $\log {\rm dist}(\cdot,I_f) \in L^1(\mu_f)$.
\end{pro}

\begin{proof}
We first characterize the [BD] condition.  Afterward, we will show how to 
`localize' it to individual points in $I_f$.  Recall from [BD] that, under the
assumption \eqref{minimality}, 
$\log \text{dist}(\cdot,I_f)$ is comparable to the function $\g^+$.  Hence it 
suffices to analyze the condition $\g^+ \in L^1(\mu_f)$.  

Suppose first that $g^-$ is finite on $I_f$.  It follows from 
Corollary 4.8 in [BD] that $g^+\in L^1(\mu_f)$, and from the bound
$g^+ \leq \g^+ \leq 0$ that $\g^+\in L^1(\mu_f)$.

Assume now that $\g^+ \in L^1(\mu_f)$.  It follows from Stokes theorem that
$$
\int_X (-\g^+) d\mu_f=O(1)+\int_X (-g^+) dd^c \g^+ \wedge T^-
\geq O(1)+\int_X d \g^+ \wedge d^c \g^+ \wedge T^-,
$$
where the $O(1)$ terms account for the fact that $g^+$ and $\g^+$ are 
plurishubharmonic only up to the addition of a smooth function.
Therefore $\g^+$ has finite energy with respect to the invariant current $T^-$.
By Stokes theorem again,
$$
\int d \g^+ \wedge d^c \g^+ \wedge T^-=O(1)+\int d \g^+ \wedge d^c \g^+ \wedge dd^c g^-
=O(1)+\int (-g^-) (dd^c \g^+)^2.
$$
Now
$$
(dd^c \g^+)^2=\sum_{p \in I_f} c_p \d_p+\l^{-2}f^* (\o^+ \wedge \o^+)
-2 \l^{-1} f^* \o^+ \wedge \o^+ +\o^+ \wedge \o^+,
$$
where $c_p>0$ and $\d_p$ denotes the Dirac mass at point $p$. Therefore
$$
\int_X d\g^+ \wedge d^c \g^+ \wedge T^-=O(1)-\sum_{p \in I_f} c_p g^-(p),
$$
so $\g^+ \in L^1(\mu_f)$ implies that $g^-$ is finite at every point of $I_f$.

We can now localize previous reasoning in the following way.
Fix $p \in I_f$ and $\chi \geq 0$ a test function supported near $p$ such that
$\chi \equiv 1$ in some small neighborhood of $p$. Thus $\f^+:=\chi \g^+$ is
comparable to $\log {\rm dist}(\cdot,p)$. Observe that 
$dd^c \f^+$ equals $\chi dd^c \g^+$ up to a smooth form.
By using Stokes theorem as above we thus get
$$
\int (-\f^+) d\mu_f=O(1)+\int (-g^-) T^+ \wedge \chi dd^c \g^+
=O(1)+\l^{-1} \int (-g^-)T^+\wedge \chi f^* \o^+.
$$
Now $\chi T^+ \wedge f^* \o^+=c \d_p+\chi f^* (T^+ \wedge \o^+)$ 
for some $c>0$.  Therefore
$$
\int (-\f^+) d\mu_f=O(1)-c\frac{g^-(p)}{\l} +
\int (-g^- \circ f^{-1})\, T^+ \wedge \o^+.
$$
The last integral is finite since $0\geq g^-\circ f^{-1} \geq \l g^-$.
Thus $\f^+ \simeq \log {\rm dist}(\cdot,p)$ is integrable with respect to 
$\mu_f$ 
if and only if 
$g^-$ is finite at point $p$.
\end{proof}

\begin{thm}
The measure $\mu_f$ is mixing.
\end{thm}

\begin{proof}
Let $\chi,\p$ be test functions. We have to show that
$$
\int_X \p (\chi \circ f^n) \,d\mu_f \longrightarrow c_{\chi} c_{\p},
\; \text{ where }
c_{\chi}=\int_X \chi \,d\mu_f \text{ and } c_{\p}=\int_X \p \,d\mu_f.
$$
It follows from the extremality of $T^+$ that the currents
$(\chi \circ f^n) T^+$ converge weakly towards $c_{\chi} T^+$.
For fixed $j$, the forms $\p \l^{-j}(f^{-j})^* \o^-$ are smooth
off the finite set $I_{f^{-j}}$. Since $T^+ \wedge \l^{-j}(f^{-j})^* \o^-$
does not charge this set, we infer
$$
\langle (\chi \circ f^n)\, T^+, \p \l^{-j}(f^{-j})^* \o^- \rangle
\longrightarrow
\langle c_{\chi} T^+, \p \l^{-j}(f^{-j})^* \o^- \rangle,
$$
as $n \rightarrow +\infty$.
Since $\l^{-j}(f^{-j})^* \o^-=\o^-+dd^c g_j^-$ with $g_j^-$ decreasing towards $g^-$,
it follows from the Monotone convergence theorem that
$$
\langle c_{\chi} T^+, \p \l^{-j}(f^{-j})^* \o^- \rangle
\longrightarrow 
\langle c_{\chi} T^+, \p T^- \rangle=c_{\chi} c_{\p},
$$
when $j \rightarrow +\infty$. Therefore it suffices to show that the
difference
$$
\langle \chi \circ f^n T^+, \p T^- \rangle
-\langle \chi \circ f^n T^+, \p \l^{-j}(f^{-j})^* \o^- \rangle
=\langle \chi \circ f^n T^+, \p dd^c (g^- -g_j^-)\rangle
$$
converges to $0$ as $j \rightarrow +\infty$ uniformly with respect to $n$.
By Stokes theorem, it suffices to uniformly control the following four 
quantities:
$$
A_{n,j}:=\langle \chi \circ f^n T^+ \wedge dd^c \p,(g^- -g_j^-) \rangle;
\;
B_{n,j}:=\langle \p dd^c (\chi \circ f^n) \wedge  T^+ ,(g^- -g_j^-) \rangle;
$$
and
$$
C_{n,j}:=\langle d \p \wedge d^c(\chi \circ f^n) \wedge  T^+,(g^- -g_j^-) \rangle;
\;
D_{n,j}:=\langle d^c\p \wedge d(\chi \circ f^n) \wedge  T^+,(g^- -g_j^-) \rangle.
$$
Observe first that
$$
|A_{n,j}| \leq ||\p||_{{\mathcal C}^2} ||\chi||_{{\mathcal C}^0} \; 
\langle \Omega \wedge T^+, (g_j^- -g^-) \rangle
$$
To control $B_{n,j}$, we observe that $dd^c (\chi \circ f^n\,T^+)$
does not charge curves.   Hence changing variables gives
$$
|B_{n,j}|=\left| \langle dd^c \chi \wedge T^+ ,\p \circ f^{-n} (g^- -g_{j+n}^-) \rangle \right|
\leq 
||\p||_{{\mathcal C}^0} ||\chi||_{{\mathcal C}^2} \; 
\langle \Omega \wedge T^+, (g_{j+n}^- -g^-) \rangle.
$$
To control $C_{n,j}$, we use the Cauchy-Schwarz inequality and obtain
\begin{eqnarray*}
|C_{n,j}| &\leq& \langle d \p \wedge d^c \p \wedge  T^+,(g^- -g_j^-) \rangle^{1/2}
\cdot 
\langle d \chi \circ f^n \wedge d^c \chi \circ f^n \wedge  T^+,(g^- -g_j^-) \rangle^{1/2} \\
&\leq& 
||\p||_{{\mathcal C}^1} ||\chi||_{{\mathcal C}^1} \; 
\langle \Omega \wedge T^+, (g_j^- -g^-) \rangle^{1/2} \cdot
\langle \Omega \wedge T^+, (g_{j+n}^- -g^-) \rangle^{1/2}.
\end{eqnarray*}
The estimation for $D_{n,j}$ is similar. This shows that $\mu_f$ is mixing.
\end{proof}

Thanks to the preceding results, we can show that $\mu_f$ is not too
concentrated.

\begin{cor}
The measure $\mu_f$ does not charge compact complex curves.
\end{cor} 

\begin{proof}
Suppose first that $\mu_f$ charges some point $p\in X$.  By Theorems 5.3 and
5.5, $p$ is a fixed point not in $I_f$ or $I_{f^{-1}}$, and $\mu_f$ is
concentrated entirely at $p$.  In particular,
both functions $\log\dist(\cdot,I_f)$ and $\log\dist(\cdot,I_{f^{-1}})$ are 
$\mu_f$-integrable.  From Proposition \ref{finiteg} we obtain that $g^+$ is
finite on $I_{f^{-1}}$ and $g^-$ is finite on $I_f$---i.e. that the [BD]
condition holds.  Theorem 4.10 from [BD] then implies that $\mu_f$ does not
charge points, which is a contradiction.

Now suppose that $\mu_f$ charges some irreducible curve $V\subset X$.  Then 
invariance of $\mu_f$ implies that $V$ cannot be critical for $f$, because 
$f(V-I_f)\subset I_{f^{-1}}$.  Invariance also implies that $\mu_f$ almost 
every point is non-wandering.  Hence 
$f^k$ restricts to an automorphism of $V$ for some $k\geq 0$.  However, the 
only mixing invariant measures for automorphisms of curves are point masses 
concentrated at fixed points, and we have already ruled out the possibility
that $\mu_f$ charges points.
\end{proof}

In section 6 we will see examples where
$g^+ \in L^1(T^- \wedge \Omega)$ but the conditions in both Theorem 4.1 and
in [BD] fail.  Observe that when the condition in [BD] is not satisfied, it is 
unclear
whether $f$ has well defined Lyapunov exponents (see proposition 5.4).

\section{Examples}

In this section, we present two examples that complement the theorems above.
The first shows that the indeterminacy orbit of a rational map can be dense
in the host manifold.  The second, which occupies the majority of the section,
shows that the invariant measure $\mu_f = T^+\wedge T^-$ can exist for a birational 
surface map even when the map fails to satisfy the equivalent conditions 
in Theorem 4.1.

\subsection{A rational map with dense indeterminacy orbits}
Let $Y=E \times E$ be a complex torus, where
$E=\C/\Z[\zeta]$ is the elliptic curve associated to a primitive root of unity
$\zeta$ of order $3,4$ or $6$. The matrix
$$
A=\left[
\begin{array}{cc} d & 1 \\ 1 & d \end{array} \right],
\; d \geq 3,
$$
preserves the lattice $\Lambda=\Z[\zeta] \times \Z[\zeta]$ and 
thus induces a holomorphic
endomorphism $g:Y \rightarrow Y$ such that
$$
\l_2(g)=(d^2-1)^2 \;  \text{ et } \;
\l_1(g)=(d+1)^2<\l_2(g).
$$
Let $\sigma:Y \rightarrow Y$ be multiplication by $\zeta$, and
let $f:X \rightarrow X$ denote the endomorphism induced by $g$ 
on the rational surface $X$ obtained
by desingularizing
the quotient $Y/\langle \sigma\rangle $, i.e. by blowing up
at fixed points of $\sigma$.

Let $a$ be such a fixed point. Since
$g$ has topological degree $\l_2(g) \geq 2$,
$g^{-1}(a)$ contains preimages other than the fixed points of $\sigma$. 
Each
point in $g^{-1}(a) \setminus Fix(\sigma)$ corresponds, in $X$, 
to a point of indeterminacy of $f$. 
Since the Lebesgue measure $\nu_Y$ of the torus $Y$ is $g$-mixing,
the preimages $(g^{-n}(a))_{n \in \N}$ are equidistributed 
with respect to $\nu_Y$ and therefore dense in $Y$.
It follows that the set
$$
I_f^{\infty}:=\bigcup_{n \in \N} f^{-n} I_f
\text{ is dense in } X.
$$
Observe also that $f$ is 1-stable: since $g$ does not  contract any curve,
neither does $f$.

\subsection{A birational surface map with constrained indeterminacy orbits}
Our second example is a variation on one due to Favre [F].
For parameters $a,b,c \in \C^*$, we consider $f=f_{abc}:\P^2 \rightarrow
\P^2$ defined by
$$
f[x:y:z]=[bcx(-cx+acy+z):
acy(x-ay+abz):abz(bcx+y-bz)].
$$
The following facts can be verified by straightforward computation.
\begin{itemize}
\item $f_{abc}$ is birational with inverse $f^{-1} = f_{a^{-1}b^{-1}c^{-1}}$.
\item $I_f = \{[a:1:0],[0:b:1],[1:0:c]\}$.
\item $f$ preserves each of the lines $\{x=0\}$, $\{y=0\}$, $\{z=0\}$
  according to the formulas
$$
[x:1:0] \mapsto [-\frac{bc}{a}x:1:0],\quad
[0:y:1]\mapsto \left[0:-\frac{ac}{b}y:1\right],\quad
[1:0:z]\mapsto \left[1:0:-\frac{ba}{c}z\right]
$$
\end{itemize}
In particular, we have $I^\infty_f, I^\infty_{f^{-1}} \subset \{xyz=0\}$ for
all $a,b,c\in\C^*$.  Let $\Omega$ denote the Fubini
Study K\"ahler form on $\P^2$.  We will spend the rest of this section
proving

\begin{thm} 
\label{breg}
Given $s>1$ and an irrational number $\theta\in\R$, let 
$f:\P^2\to\P^2$ be the birational map $f=f_{abc}$ with $a=i$, 
$b=-se^{2\pi i\theta}$, $c=i/s$.  Then $f$ is $1$-stable with invariant
currents $T^\pm = \Omega + dd^c g^\pm$ such that 
$g^+ \in L^1(\Omega\wedge T^-)$.  Moreover, for suitable choices of the
irrational number $\theta$, the condition 3 in Theorem 4.1 fails.
\end{thm}

In other words, the measure $\mu = T^+\wedge T^-$ is well-defined even though
the conditions in Theorem 4.1 sometimes fail.  

\proof
To see that $f$ is $1$-stable, observe that the hypothesis $s>1$ implies
that $f^{-n}I_f\cap\{x=0\} \subset \{[0:y:1]:|y|>1\}$ for all $n\geq 0$, whereas
$I_{f^{-1}}\cap \{x=0\} \subset \{[0:y:1]:|y|<1\}$.  In particular
$f^{-n}I_f \cap I_{f^{-1}}\cap \{x=0\} = \emptyset$.  Similarly, 
$f^{-n}I_f \cap I_{f^{-1}}\cap \{y=0\} = \emptyset$.  Finally, for each $n\geq
0$, we have $f^{-n} I_f\cap \{z=0\} = [ie^{-2\pi i n\theta}:1:0]$, and since
$\theta$ is irrational, these points never coincide with
$I_{f^{-1}}\cap\{z=0\} = [-i:1:0]$.  We conclude that 
$f^{-n} I_f\cap I_{f^{-1}}=\emptyset$ for all $n\geq 0$.  That is,
$f$ is $1$ stable.  In particular $\lambda := \r_{f^*} = 2$.

To prove that $g^+\in L^1(\Omega\wedge T^-)$, we apply Proposition 5.2 with
$S$ equal to the current of integration $[x=0]$ over the line $\{x=0\}$.
Observe that $\Omega\wedge [x=0]$ is just area measure on $\{x=0\}$.  
Because $g^-$ is qpsh it follows that either $g^-$ is integrable with respect to
$\Omega\wedge[x=0]$, or $g^-|_{\{x=0\}} \equiv - \infty.$  The latter
is far from true, however.  One can compute directly, for instance,
that $g^->-\infty$ at the fixed point $[0:0:1]$.  Thus 
$g^-\in L^1(\Omega\wedge [x=0])$.

It follows from standard arguments that local potentials for $T^-$ must
be harmonic on any open set $U\subset\{x=0\}$ such that 
\begin{itemize}
\item iterates of $f^{-1}$ form a normal family on $U$, and
\item $U\cap f^n I_{f^{-1}} = \emptyset$ for all $n\geq 0$.
\end{itemize}
The only point in $\{x=0\}$ where iterates of $f^{-1}$ fail to act 
normally is the fixed point $[0:0:1]$.  Hence
$$
\supp\, ([x=0]\cap T^-) \subset \overline{\cup_{n\geq 0} f^n I_{f^{-1}}}
$$
is a compact subset of $\{[0:y:1]:|y|<1\}$.  Replacing $f^{-1}$ with $f$, 
the same reasoning shows that local potentials for $T^+$ 
are harmonic on $\{[0:y:1]:|y|<1\}$.  Thus $g^+$ is uniformly bounded
on $\supp\,([x=0]\cap T^-)$, and it follows that $g^+\in L^1([x=0]\wedge T^-)$.
Therefore by Proposition 5.2, $g^+\in L^1(\Omega\wedge T^-)$.

To see that condition 3 in Theorem 4.1 fails for suitably chosen $\theta$,
let $h:\R^+\to\R^+$ be a function decreasing rapidly to 0.  By a Baire category
argument, one can find irrational $\theta$ such that 
$$
2n_j\theta \mod 1 < h(n_j)
$$ 
for infinitely many $n_j\in\N$.  Thus, if we set $p^+ = [i:1:0] \in I_f$
and $p^- = [-i:1:0]\in I_f$, we obtain
\begin{eqnarray*}
\sum_{n=0}^\infty \frac{1}{\lambda^{2n}} \log\dist(f^{-n}I_f,I_{f^{-1}})
& \leq & 
\sum_{n=0}^\infty \frac{\log\dist(f^n(p^-),p^+)}{2^{2n}} \\
& \leq &  
C\sum_{n=0}^\infty \frac{\log|e^{-2\pi n i}+1|}{2^{2n}} \leq
C\sum_{j=0}^\infty \frac{\log h(n_j)}{2^{2n_j}}.
\end{eqnarray*}
The last sum diverges to $-\infty$ if we take e.g. $h(x) = 2^{-2^{2n}}$,
and condition 3 in Theorem 4.1 then fails.
\qed

\vskip 1cm

Jeffrey DILLER

Department of Mathematics

University of Notre Dame

Notre Dame, Indiana 46656

USA

diller.1@nd.edu

\vskip1cm

Vincent Guedj

Universit\'e Aix-Marseille 1, LATP

13453 MARSEILLE Cedex 13

FRANCE

guedj@cmi.univ-mrs.fr

\end{document}